\documentclass[a4paper,12pt]{article}
\usepackage{amsmath}
\usepackage{amsthm}
\usepackage{amsfonts}
\usepackage{mathtools}
\usepackage{amssymb}
\usepackage{ulem}
\usepackage{color}
\usepackage{graphicx} 
\usepackage{tikz}
\usepackage{float}
\usetikzlibrary{intersections, calc, arrows.meta}
\usepackage{mleftright}
\usepackage{braket}
\usepackage[margin=25mm]{geometry}

\theoremstyle{definition}
\newtheorem{prop}{Proposition}[section] 

\newtheorem*{rem}{Remark}               

\newtheorem{thm}[prop]{Theorem}
\newtheorem{lem}[prop]{Lemma}

\numberwithin{equation}{section}

\everymath{\displaystyle}

\title{Asymptotic Formula of the Lifespan of Classical Solutions to 1D Quasilinear Wave Equations}
\date{}
\author{Yuusuke Sugiyama and Taro Yamanoi}

\begin{document}

\maketitle

\begin{abstract}
	We study the lifespan \(T^*(\varepsilon)\) of classical solutions to the Cauchy problem
	for the one-dimensional quasilinear wave equation
	$u_{tt}-(c(u)^2u_x)_x=0$,
	with initial data of size $\varepsilon>0$.
	We assume that $c'(\theta)\sim a|\theta|^{p-2}$ as $\theta \to 0$,
	where \(a>0\) and \(p>3\) or $p=2$.
	We obtain explicit asymptotic formulas for
	$T^*(\varepsilon)$ as $\varepsilon\downarrow0$.
	We treat both the case $\psi \in L^1(\mathbb R)$
	and the case
	$\psi=\psi_1+\Braket{x}^{-\alpha}$,
	where \(\psi_1\in L^1(\mathbb R)\) and \(0<\alpha\le1\).
	Our results identify the dependence of the lifespan scale on the decay exponent,
	including the borderline case $\alpha=1$.
	The proof combines Riccati equations along characteristics
	with a scaling argument inspired by Lindblad.
\end{abstract}

\section{Introduction}

\subsection{Problem}

In this paper, we consider the following Cauchy problem of the one-dimensional quasilinear wave equation.
\begin{equation} \label{Eq:Main}
	\begin{cases}
		u_{tt} - (c(u)^2u_x)_x = 0, \hspace{8pt}  (t,x) \in [0,T]\times \mathbb{R}, \\
		u(0,x)=\varepsilon \varphi(x),\,u_t(0,x)=\varepsilon \psi(x), \hspace{8pt}  x \in \mathbb{R}.
	\end{cases}
\end{equation}
Let $p > 3$ or $p=2$, and let $U\subset \mathbb{R}$ be a neighborhood of the origin.
Throughout this paper, we assume that
\begin{equation*}
	c(0)=1 \hspace{10pt} \text{and} \hspace{10pt} c(\theta)>0 \hspace{10pt} \text{for all }\theta \in U.
\end{equation*}
In addition, we assume that
\begin{equation}\label{Eq:unif cvg}
	c \in C^2(U)  \hspace{10pt} \text{and} \hspace{10pt} \sup_{|\theta| \leq M}(|\xi^{-(p-2)}c'(\xi \theta)-a|\theta|^{p-2}|) \rightarrow 0 \hspace{8pt} \text{as }\xi\downarrow0
\end{equation}
for all $M>0$, where $a>0$ is a positive constant.
Examples of functions $c(\theta)$ satisfying the assumptions include the following:
\begin{align*}
	 & c(\theta)=(1+A|\theta|^{p-2}\theta)^{\frac{1}{2}},       \\
	 & c(\theta)=e^{A|\theta|^{p-2}\theta},                     \\
	 & c(\theta)=1+A|\theta|^{p-2}\theta+B|\theta|^{q-2}\theta,
\end{align*}
where \(A,B>0\), \(q>p\), and \(q>3\).
The existence and uniqueness of a local-in-time solution to the Cauchy problem \eqref{Eq:Main} is known for $\varphi \in C_b^2$ and $\psi \in C_b^1$ (e.g. Hartman and Wintner \cite{HW} and see Remark).
The blow-up alternative for equation \eqref{Eq:Main} implies that
a solution can be extended as long as $\|u_t(t)\|_{L^\infty}+\|u_x(t)\|_{L^\infty}$
remains bounded and the range of $u$ stays in a compact subset of $U$.
From this blow-up alternative, we define the life-span of solutions for the Cauchy problem (1.1) by
\begin{align*}
	T^*(\varepsilon) = \sup \biggl\{ T>0 \;\biggm|\; \sup_{t\in [0,T)} \Bigl(
	 & \|u_t(t)\|_{L^\infty}+\|u_x(t)\|_{L^\infty} \Bigr) < \infty \hspace{10pt} \text{and} \hspace{10pt} \overline{u([0,T)\times \mathbb{R})} \subset U\biggr\}.
\end{align*}

We now state the main results of this paper. The first theorem gives an asymptotic formula for the lifespan when the initial datum $\psi$ is integrable.
\begin{thm}
	Let $\varphi \in C_b^2, \psi \in C_b^1 \cap L^1$ and further suppose that $\varphi, \varphi_x, \psi \rightarrow 0$ as $|x| \rightarrow \infty.$
	Moreover, we assume that there exists a point $x_0 \in \mathbb{R}$ such that one of the following two conditions is satisfied:
	\begin{align*}
		(\mathrm{i})  & \,\,\psi(x_0)+\varphi_x(x_0)>0, \\[7pt]
		(\mathrm{ii}) & \,\,\psi(x_0)-\varphi_x(x_0)>0.
	\end{align*}
	Then
	\begin{equation*}
		\lim_{\varepsilon \rightarrow 0}T^*(\varepsilon) \varepsilon^{p-1}=\Lambda_1^{-1},
	\end{equation*}
	where
	\begin{align*}
		                                     & \Lambda_1=\frac{a}{2^{p-1}} \max\{ \sup_{x \in \mathbb{R}}(B_{+}(x)|A_{+}(x)|^{p-2}), \sup_{x \in \mathbb{R}}(B_{-}(x)|A_{-}(x)|^{p-2})\}, \\[9pt]
		A_{+}(x)=\int_{-\infty}^{x}\psi(y)dy & + \varphi(x), \hspace{5pt}  A_{-}(x)=\int_{x}^{\infty}\psi(y)dy + \varphi(x),  \hspace{5pt} B_{\pm}(x)= \psi(x) \pm \varphi_x(x).
	\end{align*}
\end{thm}
Our second main theorem concerns nonintegrable initial
data $\psi$ with a prescribed algebraic rate of spatial decay.
We derive asymptotic formulas for the lifespan,
with a scaling determined by the spatial decay exponent.

\begin{thm}
	Let $\varphi \in C_b^2, \psi=\psi_1+\psi_2,$ where $\psi_1 \in C_b^1 \cap L^1$ and $\psi_2=\Braket{x}^{-\alpha}$
	with $0 < \alpha \leq 1.$
	Moreover, we suppose that $\varphi, \varphi_x, \psi \rightarrow 0$ as $|x| \rightarrow \infty.$
	Furthermore, we assume that there exists a point $x_0 \in \mathbb{R}$ such that one of the two conditions $(\mathrm{i}), (\mathrm{ii})$ in Theorem 1.1 is satisfied.
	Then
	\begin{align*}
		\lim_{\varepsilon \rightarrow 0}T^*(\varepsilon)\varepsilon^{q}                          & =\Lambda_2^{-1} \hspace{10pt}  \text{if} \hspace{10pt}0<\alpha<1, \\
		\lim_{\varepsilon \rightarrow 0}T^*(\varepsilon)\varepsilon^{p-1}|\log\varepsilon|^{p-2} & =\Lambda_3^{-1} \hspace{10pt} \text{if}  \hspace{10pt}\alpha=1,
	\end{align*}
	where $q=\frac{p-1}{(1-\alpha)(p-2)+1}$ and
	\begin{align*}
		\Lambda_2 & =\left(\frac{2^{-\alpha(p-2)}a}{2k(1-\alpha)^{p-2}} \max\left({\sup\limits_{x \in \mathbb{R}} B_+(x),\sup\limits_{x \in \mathbb{R}} B_-(x)}\right)\right)^{\frac{1}{k}}, \\[4pt]
		\Lambda_3 & =\frac{a(p-1)^{p-2}}{2^{p-1}}\max\left(\sup_{x \in \mathbb{R}}B_+(x), \sup_{x \in \mathbb{R}}B_-(x)\right),
	\end{align*}
	where $k=(1-\alpha)(p-2)+1.$
\end{thm}
These limit formulas imply upper and lower bounds for
\(T^*(\varepsilon)\) with the optimal dependence on
\(\varepsilon\). Even these bounds are new.

\subsection*{Notation}
$C_b^m$ is the set of $m$-times differentiable functions such that the functions themselves and all their derivatives up to order $m$ are bounded and continuous on $\mathbb{R}$.
Also, throughout this paper, $C$ and $C_j$ denote generic positive constants whose values may change from line to line.

\subsection{Known results}
Lax \cite{Lax64} and John \cite{John} investigated more general $2\times 2$ and $n \times n$ systems of hyperbolic conservation laws, respectively.
By applying these results, we can obtain
\begin{equation*}
	C\varepsilon^{-1} \leq T^*(\varepsilon) \leq C'\varepsilon^{-1}
\end{equation*}
for the case $c(u)>0\,,\,c'(0) \not = 0$ (that is, $p=2$).
This case is referred  to as genuine nonlinearity, which was introduced by Lax \cite{Lax57}.
Furthermore, to obtain more detailed information on $T^*(\varepsilon)$, H\"{o}rmander \cite{Hormander} derived the limit
\begin{equation*}
	\lim_{\varepsilon \to 0} T^*(\varepsilon) \varepsilon = \Lambda^{-1}.
\end{equation*}
In his work, initial data with compact support are considered.
The case $c'(0)=0$, which is close to  the linear wave $u_{tt}-c^2u_{xx}=0$, has also been widely investigated.
For the case $u_{tt}-c(u_x)^2u_{xx}=0$, where $c(\theta) = (1 + A|\theta|^{p-2}\theta)^{1/2}$ with $A > 0$ and $p > 1$,
Haruyama and Takamura \cite{Haruyama} established the upper bound
\begin{equation*}
	T^*(\varepsilon) \leq C \varepsilon^{-(p-1)}
\end{equation*}
by using the method of successive approximations.
They proved this under the assumption that the initial data are non-negative and have compact support.
However, the existence of local-in-time solutions is not known for $1 < p < 2$.
Moreover, for more general $n \times n$ hyperbolic systems, Li, Zhou and Kong \cite{LZK} established
\begin{equation}\label{Eq:T est}
	C \varepsilon^{-(p-1)} \le T^*(\varepsilon) \le C' \varepsilon^{-(p-1)}
\end{equation} for integer $p \ge 2$,
imposing a decay of $|x|^{-(1+\mu)}$ with $\mu > 0$ on the initial data.
Also, Sugiyama and Yamanoi \cite{S.Y} established estimate \eqref{Eq:T est} for the equation $u_{tt} - c(u_x)^2u_{xx} = 0$,
removing the decay rate assumption and the restriction that $p \geq 2$ is an integer.
As in the case where $c'(0) \neq 0$,
the limit
\begin{equation*}
	\lim_{\varepsilon \to 0} T^*(\varepsilon) \varepsilon^{p-1} = \Lambda^{-1}
\end{equation*}
has been derived by Li and Kong \cite{LK} and Kong \cite{Kong}.
They also considered more general $n \times n$ hyperbolic systems,
imposing a decay of $|x|^{-(1+\mu)}$ with $ \mu>0$ and $|x|^{-1}$ on the initial data, respectively, while restricting $p \geq 2$ to integers.
Let $\Psi$ be a primitive of $\psi$.
If $\varphi\in C_b^2$, $\psi\in C_b^1\cap L^1$,
and $\varphi(x), \varphi_x(x), \psi(x), \Psi(x)=O(|x|^{-1})$
as $|x|\to\infty$, their results can be applied to the Cauchy problem \eqref{Eq:Main}.
In particular, the decay condition on the primitive \(\Psi\)
requires that $\int_{\mathbb R}\psi(x)\,dx=0$.
Theorem 1.1 of this paper constitutes a partial improvement over their results,
in the sense that both the decay rate assumption and the integer-order restriction on $p$ are removed.
Their proof is based on a priori estimates for weighted $L^\infty$ norms via a bootstrap argument.
In this paper, we give the proof by using the representation of solutions along characteristics obtained via Riemann invariants,
together with the scaling transformation by Lindblad \cite{Lindblad}.

For semilinear wave equations, the dependence of the lifespan
rate on the spatial decay of the initial data, as in
Theorem 1.2, was already known; see, e.g., Asakura \cite{Asakura}, Tsutaya \cite{Tsutaya} and Takamura \cite{Takamura}.
The first limit formula for the lifespan of solutions to
semilinear wave equations was obtained by Lindblad \cite{Lindblad}.
His result assumes that the initial data have compact
support and hence are integrable.
However, when nonintegrable initial data lead to a different
lifespan rate, no corresponding limit formula is known.

For quasilinear wave equations, even the upper and lower
bounds for $T^*(\varepsilon)$ obtained in Theorem 1.2 are new.
Moreover, Theorem 1.2 gives the corresponding limit formulas.

For semilinear heat equations, Lee and Ni \cite{LN}
obtained lifespan estimates whose rates depend on
the spatial decay of the initial data.
Dickstein \cite{Dickstein} further obtained corresponding
limit formulas for the lifespan as the size of
the initial data tends to zero.

\begin{rem}[Existence for $C^1$ solutions with $2< p \leq 3$]
	Suppose that \(2<p\le3\), \(c\in C^1(U)\), and the
	asymptotic condition in \eqref{Eq:unif cvg} holds.
	If \(\psi\in L^1\), introducing a primitive \(v\) of
	\(u_t\) gives the following \(2\times2\) system:
	\begin{equation}\label{Eq:integrated system}
		\begin{cases}
			u_t-v_x=0,                    \\
			v_t-c(u)^2u_x=0,              \\
			u(0,x)=\varepsilon\varphi(x), \\
			v(0,x)=\displaystyle\varepsilon
			\int_{-\infty}^x\psi(y)\,dy.
		\end{cases}
	\end{equation}
	This is a strictly hyperbolic system with \(C^1\)
	coefficients, and hence has a unique local-in-time
	\(C^1\) solution; see \cite{HW}.
	The conclusion of Theorem 1.1 remains valid for the
	lifespan of this solution.

	When \(\psi\notin L^1\), we instead consider the following
	system for the Riemann invariants \(r,s\) defined in
	\eqref{Eq:Riemann invariants}:
	\begin{equation}\label{Eq:closed Riemann system}
		\begin{cases}
			r_t-
			c\!\left(G^{-1}\!\left(\dfrac{r-s}{2}\right)\right)r_x=0, \\[8pt]
			s_t+
			c\!\left(G^{-1}\!\left(\dfrac{r-s}{2}\right)\right)s_x=0, \\[8pt]
			u=G^{-1}\!\left(\dfrac{r-s}{2}\right),
		\end{cases}
	\end{equation}
	where \(G(u)=\int_0^u c(\theta)\,d\theta\).
	The initial data are
	\begin{equation*}
		r(0,x)=\varepsilon\int_0^x\psi(y)\,dy
		+G(\varepsilon\varphi(x)),
		\qquad
		s(0,x)=\varepsilon\int_0^x\psi(y)\,dy
		-G(\varepsilon\varphi(x)).
	\end{equation*}
	To construct a local solution, we set
	\begin{equation*}
		R(t,x)=r(t,x)-\varepsilon\int_0^x\psi(y)\,dy,
		\qquad
		S(t,x)=s(t,x)-\varepsilon\int_0^x\psi(y)\,dy.
	\end{equation*}
	Then
	\begin{equation*}
		\begin{cases}
			R_t-c(u)R_x=\varepsilon c(u)\psi(x),  \\
			S_t+c(u)S_x=-\varepsilon c(u)\psi(x), \\[4pt]
			u=G^{-1}\!\left(\dfrac{R-S}{2}\right),
		\end{cases}
	\end{equation*}
	with
	\begin{equation*}
		R(0,x)=G(\varepsilon\varphi(x)),
		\qquad
		S(0,x)=-G(\varepsilon\varphi(x)).
	\end{equation*}
	This is a strictly hyperbolic system with \(C^1\)
	coefficients and initial data in \(C_b^1\).
	The standard method of characteristics therefore
	gives a unique local-in-time \(C^1\) solution.
	The two systems have the same maximal existence time,
	which agrees with \(T^*(\varepsilon)\) defined above
	for the reconstructed \(u\), with the same continuation
	criterion.
	Under the other assumptions of Theorem 1.2,
	the same lifespan formulas hold for \(C^1\) solutions
	of \eqref{Eq:closed Riemann system} when \(2<p\le3\).

	In these formulations, \(u\) satisfies
	\eqref{Eq:Main} in the distributional sense, and
	\(T^*(\varepsilon)\) denotes the lifespan of the
	\(C^1\) solution to the corresponding first-order
	system, with the same continuation criterion
	as above.
\end{rem}

\section{Preliminaries}
We define the Riemann invariants as follows:
\begin{equation}\label{Eq:Riemann invariants}
	\begin{cases}
		r(t,x)=\int_{0}^{x}u_t(t,y)dy + \int_{0}^{t}c(u(\tau,0))^2u_x(\tau,0)d\tau + G(u), \\[10pt]
		s(t,x)=\int_{0}^{x}u_t(t,y)dy + \int_{0}^{t}c(u(\tau,0))^2u_x(\tau,0)d\tau - G(u),
	\end{cases}
\end{equation}
where $G(u)=\int_{0}^{u}c(\theta)d\theta.$
Then
\begin{equation}\label{Eq:invariance}
	\begin{cases}
		r_t-c(u)r_x=0, \\
		s_t+c(u)s_x=0.
	\end{cases}
\end{equation}
For $c(u)$, the positive and negative characteristics are solutions to
the following differential equations:
\begin{equation*}
	\begin{cases}
		\frac{dx_\pm}{dt}(t)=\pm c(u(t,x_\pm(t))), \\[10pt]
		x_\pm(\tau)=y.
	\end{cases}
\end{equation*}
We denote the characteristics through $(\tau,y)$ by $x_\pm(t;\tau,y)$,
and we simply write $x_\pm(t)$ for $x_\pm(t;\tau,y)$ unless there is a risk of confusion.
Using the chain rule, by \eqref{Eq:invariance}, we obtain
\begin{equation}\label{Eq:Invariance1}
	\begin{cases}
		\frac{d}{dt}r(t,x_{-}(t))=0, \\[10pt]
		\frac{d}{dt}s(t,x_{+}(t))=0.
	\end{cases}
\end{equation}
Integrating \eqref{Eq:Invariance1} over $[0,t]$, we have
\begin{equation}\label{Eq:Invariance2}
	\begin{cases}
		r(t,x_{-}(t))=r(0,x_{-}(0)), \\[10pt]
		s(t,x_{+}(t))=s(0,x_{+}(0)).
	\end{cases}
\end{equation}
That is, $r$ and $s$ are invariant on $x_{-}(t)$ and $x_{+}(t)$, respectively.
From the equation in \eqref{Eq:Main}, differentiating both sides of each equation in \eqref{Eq:invariance} with respect to $x$, we obtain
\begin{equation}
	\begin{cases}
		(\partial_t-c\partial_x)r_x=\frac{c'}{2c}r_x(r_x-s_x), \\[10pt]
		(\partial_t+c\partial_x)s_x=\frac{c'}{2c}s_x(s_x-r_x).
	\end{cases}
\end{equation}
Let $F_1=\sqrt{c}r_x$ and $F_2=\sqrt{c}s_x$.
Then we obtain the following Riccati-type differential equations:
\begin{align}
	 & \quad \frac{d}{dt}F_1(t,x_{-}(t))=\gamma F_1^2,  \label{Eq:Riccati1} \\[3pt]
	 & \quad \frac{d}{dt}F_2(t,x_{+}(t))=\gamma F_2^2,  \label{Eq:Riccati2}
\end{align}
where $\gamma=\gamma(u)=\frac{c'(u)}{2c(u)^{\frac{3}{2}}}$.
In this paper, we prove Theorems 1.1 and 1.2 based on these Riccati-type equations.

\section{Proof of Theorem 1.1}
In this section, we assume that the initial data satisfy the assumptions of Theorem 1.1.
Then the following conservation laws hold for the Riemann invariants $r,s$.
\begin{prop}
	For $t \in [0,T^*(\varepsilon))$,
	\begin{equation}
		\|r(t)\|_{L^\infty}=\|r(0)\|_{L^\infty}, \quad \|s(t)\|_{L^\infty}=\|s(0)\|_{L^\infty}.
	\end{equation}
\end{prop}
\begin{proof}
	This follows easily from \eqref{Eq:Invariance2} and the fact that the maps from $x$ to $x_{\pm}(t;0,x)$ for all $t \in [0,T^*(\varepsilon))$ are bijections.
\end{proof}
Furthermore, the following proposition holds.
\begin{prop}
	There exists $\varepsilon_0>0$ such that if $0 < \varepsilon \leq \varepsilon_0$, then, for $t \in [0,T^*(\varepsilon))$,
	\begin{align}
		\|u(t)\|_{L^\infty}      & \leq C_1\varepsilon,  \label{Eq:u_x}                \\
		\|c'(u(t))\|_{L^\infty}  & \leq C_1\varepsilon^{p-2}, \label{Eq:c_prime_upper} \\
		\|c(u(t))-1\|_{L^\infty} & \leq C_1\varepsilon^{p-1}. \label{Eq:c_minus_1}
	\end{align}
\end{prop}
\begin{proof}
	This proposition follows by the exact same argument as in Sugiyama and Yamanoi \cite{S.Y}.
	Regarding the estimate for $c'(u)$, note that for any $\eta>0$,
	\begin{equation*}
		|c'(u)| \leq \left(\eta+a\left|\frac{u}{\varepsilon}\right|^{p-2}\right)\varepsilon^{p-2}
	\end{equation*}
	follows from the assumption \eqref{Eq:unif cvg} provided that $\varepsilon_0$ is chosen sufficiently small.
\end{proof}
In particular, by \eqref{Eq:c_minus_1},
\begin{equation}\label{Eq: c near 1}
	\frac{1}{2} \leq c(u) \leq 2.
\end{equation}
We give the proof by estimating the limit superior and the limit inferior of $T^*(\varepsilon) \varepsilon^{p-1}$.
For this purpose, we employ the scaling technique used in Lindblad \cite{Lindblad}.
Solving the Riccati-type equation \eqref{Eq:Riccati1}, we have
\begin{align*}
	F_1(t) & =F_1(t,x_{-}(t;0,x))                                                            \\[7pt]
	       & =\frac{F_1(0)}{1-F_1(0)\int_{0}^{t}\gamma(u(\sigma,x_{-}(\sigma;0,x)))d\sigma}.
\end{align*}
Similarly,
\begin{equation*}
	F_2(t)=\frac{F_2(0)}{1-F_2(0)\int_{0}^{t}\gamma(u(\sigma,x_{+}(\sigma;0,x)))d\sigma}.
\end{equation*}
We set
\begin{align*}
	D^{(1)}(t,x)=1-F_1(0)\int_{0}^{t}\gamma(u(\sigma,x_{-}(\sigma;0,x)))d\sigma, \\
	D^{(2)}(t,x)=1-F_2(0)\int_{0}^{t}\gamma(u(\sigma,x_{+}(\sigma;0,x)))d\sigma,
\end{align*}
and
\begin{equation*}
	D^{(j)}_\varepsilon(\tau,x)=D^{(j)}(\varepsilon^{-(p-1)}\tau,x).
\end{equation*}

\begin{lem}
	We fix $R>0, \tau_0 >0$ and set $I_\varepsilon=[0,\min(\tau_0, T^*(\varepsilon)\varepsilon^{p-1})).$
	Then,
	\begin{align*}
		\sup_{(\tau,x)\in I_\varepsilon \times [-R,R]}|D^{(1)}_\varepsilon(\tau,x) - D_+(\tau,x)| \rightarrow 0, \\
		\sup_{(\tau,x)\in I_\varepsilon \times [-R,R]}|D^{(2)}_\varepsilon(\tau,x) - D_-(\tau,x)| \rightarrow 0,
	\end{align*}
	as $\varepsilon \rightarrow 0$, where
	\begin{equation*}
		D_\pm(\tau,x)=1-B_\pm(x)K_\pm(x)\tau, \hspace{10pt} K_\pm(x)=\frac{a}{2^{p-1}}|A_\pm(x)|^{p-2}.
	\end{equation*}
\end{lem}
\begin{proof}
	We show that, for any $\eta>0$, there exists $\varepsilon_0>0$ such that
	\begin{align*}
		|D^{(1)}_\varepsilon(\tau,x) - D_+(\tau,x)| < \eta, \\
		|D^{(2)}_\varepsilon(\tau,x) - D_-(\tau,x)| < \eta,
	\end{align*}
	for $(\tau,x) \in [0,\min(\tau_0, T^*(\varepsilon)\varepsilon^{p-1})) \times [-R,R].$
	We fix $\delta>0$ arbitrarily and set $I_{\delta,\varepsilon} = [\delta, \min(\tau_0, T^*(\varepsilon)\varepsilon^{p-1}))$.
	By the definition of Riemann invariants \eqref{Eq:Riemann invariants},
	\begin{equation*}
		\begin{cases}
			r(0,x)=\varepsilon \int_{0}^{x}\psi(y)dy + G(\varepsilon \varphi(x)), \\[12pt]
			s(0,x)=\varepsilon \int_{0}^{x}\psi(y)dy - G(\varepsilon \varphi(x)).
		\end{cases}
	\end{equation*}
	We assume that $(\tau, x) \in I_{\delta,\varepsilon} \times [-R,R]$ and set $t=\varepsilon^{-(p-1)}\tau.$
	We consider $x_-(\sigma)$ through $(0,x)$ and
	$x_+(\sigma)$ intersecting $x_-(\sigma)$ at time $t$, and let $(t,X)$ be their point of intersection.
	We set $x_+(0)=y_\varepsilon$, since the initial value of $x_+(\sigma)$ depends on $\varepsilon$.
	\begin{figure}[htbp]
		\centering
		\begin{tikzpicture}[scale=1.0]
			\draw[->] (-3.7,0) -- (6,0) node[right] {$x$};
			\draw[->] (0,-0.5) -- (0,4) node[above] {$\sigma$};

			\draw[thick, black, name path=minus] (3.8,0) node[below]{$x$} to[out=140, in=-80] (-0.5,4.0) node[left]{$x_-(\sigma)$};

			\draw[thick, black!60!black, name path=plus] (-1.2,0) node[below]{$y_\varepsilon= x_{+}(0)$} to[out=40, in=-115] (2.5,4.0) node[right]{$x_{+}(\sigma)$};

			\path [name intersections={of=minus and plus, by=P}];
			\fill (P) circle (2pt) node[right, xshift=10pt] {$(t, X)$};

		\end{tikzpicture}
		\caption{Characteristics}
	\end{figure}
	From $2G(u)=r-s$ and \eqref{Eq:Invariance2},
	\begin{align*}
		2G(u(t,X)) & =r(0,x)-s(0,y_\varepsilon)                                                                                        \\
		           & = \int_{y_\varepsilon}^{x} \varepsilon \psi(y)dy+G(\varepsilon \varphi(x))+G(\varepsilon \varphi(y_\varepsilon)).
	\end{align*}
	Then
	\begin{align*}
		X=x-\int_{0}^{t}c(u(\sigma,x_{-}(\sigma)))d\sigma, \\
		X=y_\varepsilon+\int_{0}^{t}c(u(\sigma,x_{+}(\sigma)))d\sigma.
	\end{align*}
	Hence, from \eqref{Eq: c near 1},
	\begin{align*}
		y_\varepsilon & \leq x-t                           \\
		              & = x-\varepsilon^{-(p-1)}\tau       \\
		              & \leq R-\varepsilon^{-(p-1)}\delta.
	\end{align*}
	Therefore, by the decay assumption on initial data, for any $\eta>0$, there exists $\varepsilon_0>0$ such that
	\begin{equation*}
		\left|\frac{G(u(t,X))}{\varepsilon}-\frac{A_+(x)}{2}\right| < \eta,
	\end{equation*}
	if $0<\varepsilon<\varepsilon_0$ uniformly with respect to $(\tau,x) \in I_{\delta,\varepsilon} \times [-R,R].$
	Moreover,
	\begin{align*}
		\left|\frac{u(t,X)}{\varepsilon}-\frac{A_+(x)}{2}\right| & \leq \frac{|u-G(u)|}{\varepsilon}+\left|\frac{G(u)}{\varepsilon}-\frac{A_+(x)}{2}\right| \\
		                                                         & \leq C\varepsilon+\eta.
	\end{align*}
	Thus,
	\begin{equation*}
		\left|\frac{u(t,X)}{\varepsilon}-\frac{A_+(x)}{2}\right| <\eta,
	\end{equation*}
	if $0<\varepsilon<\varepsilon_0$ uniformly with respect to $(\tau,x) \in I_{\delta,\varepsilon} \times [-R,R].$
	Hence, by the assumption on $c(u)$ and the continuity of the solution,
	\begin{equation*}
		|\varepsilon^{-(p-2)}\gamma(u(t,X))-K_{+}(x)| <\eta.
	\end{equation*}
	Here, for $(\tau,x) \in I_\varepsilon \times [-R,R],$
	\begin{align*}
		|D^{(1)}_\varepsilon(\tau,x)-D_+(\tau,x)| & = \left|\frac{F_1(0)}{\varepsilon}\int_{0}^{\tau}\varepsilon^{-(p-2)}\gamma(u)d\sigma-B_+(x)K_+(x)\tau\right| \\[7pt]
		                                          & \leq \left| \frac{F_1(0)}{\varepsilon} \right|\int_{0}^{\tau}|\varepsilon^{-(p-2)}\gamma(u)-K_+(x)|d\sigma    \\
		                                          & +\left| \frac{F_1(0)}{\varepsilon}-B_+(x) \right| K_+(x)\tau.
	\end{align*}
	For any \(\eta>0\), by choosing \(\delta>0\)
	sufficiently small, we have
	\[
		\int_0^\delta
		\left|\varepsilon^{-(p-2)}\gamma(u)-K_+(x)\right|
		\,d\sigma<\eta
	\]
	by the uniform boundedness of
	\(\varepsilon^{-(p-2)}\gamma(u)\) and \(K_+(x)\).
	For this fixed \(\delta\),
	we also have
	\[
		\int_\delta^\tau
		\left|\varepsilon^{-(p-2)}\gamma(u)-K_+(x)\right|
		\,d\sigma<\eta
	\]
	by the uniform convergence above,
	provided that \(\varepsilon_0>0\) is sufficiently small.
	By the boundedness of the initial data and the assumption on $c(\theta)$,
	it is clear that there exist $M>0$ and $\varepsilon_0>0$ such that
	\begin{align*}
		 & \left|\frac{F_1(0)}{\varepsilon}\right|         \leq M,   \\
		 & \left|\frac{F_1(0)}{\varepsilon}-B_+(x)\right|  \leq \eta
	\end{align*}
	if $0< \varepsilon < \varepsilon_0$.
	Therefore,
	\begin{equation*}
		|D^{(1)}_\varepsilon(\tau,x)-D_+(\tau,x)| < \eta
	\end{equation*}
	if $0<\varepsilon< \varepsilon_0$ uniformly with respect to $(\tau,x) \in I_\varepsilon \times [-R, R].$
	The same holds for $D^{(2)}_\varepsilon(\tau,x).$
\end{proof}
Noting that, by the assumption on the initial data, $\frac{F_j(0)}{\varepsilon}$ is bounded for both $j = 1, 2$
and its $L^\infty$ norm is positive independently of $\varepsilon$ for at least one $j=1,2$,
it suffices to analyze $D^{(j)}_\varepsilon(\tau, x)$ to derive an estimate for $T^*(\varepsilon)$.
First, we derive an estimate for the limit superior of $T^*(\varepsilon)\varepsilon^{p-1}.$
\begin{lem}
	\begin{equation*}
		\limsup_{\varepsilon \rightarrow 0}T^*(\varepsilon)\varepsilon^{p-1} \leq \Lambda_1^{-1}.
	\end{equation*}
\end{lem}
\begin{proof}
	We present the proof under the assumption $\Lambda_1=\frac{a}{2^{p-1}}\sup_{x \in \mathbb{R}}(B_{+}(x)|A_{+}(x)|^{p-2})$.
	The case $\Lambda_1=\frac{a}{2^{p-1}} \sup_{x \in \mathbb{R}}(B_{-}(x)|A_{-}(x)|^{p-2})$ is completely analogous using $D^{(2)}_\varepsilon$.
	We fix \(0<\delta<\Lambda_1/2\) arbitrarily.
	By the definition of $\Lambda_1$, there exists $x \in \mathbb{R}$ such that
	\begin{equation*}
		B_{+}(x)K_{+}(x) \geq \Lambda_1-\delta.
	\end{equation*}
	We assume that for any $\varepsilon_0 >0$, there exists $\varepsilon>0$ such that
	\begin{equation*}
		0 <\varepsilon< \varepsilon_0 \hspace{10pt} \text{and} \hspace{10pt} T^*(\varepsilon)\varepsilon^{p-1} > \frac{1}{\Lambda_1-2\delta}.
	\end{equation*}
	Then we can choose a sequence of positive real numbers $\{\varepsilon_k\}_k$ such that
	\begin{equation*}
		\varepsilon_k \rightarrow 0 \hspace{7pt} \text{as} \hspace{7pt} k \rightarrow \infty, \hspace{10pt} \text{and} \hspace{10pt} T^*(\varepsilon_k)>\frac{\varepsilon_k^{-(p-1)}}{\Lambda_1-2\delta} \eqqcolon \varepsilon_k^{-(p-1)}\tau_1\eqqcolon t_1.
	\end{equation*}
	Consequently,
	\begin{equation}\label{Eq:lim D}
		\lim_{k \rightarrow \infty}D^{(1)}_{\varepsilon_k}(\tau_1,x) \geq 0
	\end{equation}
	must hold by $D^{(1)}_{\varepsilon_k}(0,x)=1>0$ and the continuity of the solution.
	On the other hand, we have
	\begin{align*}
		\lim_{k \rightarrow \infty} D^{(1)}_{\varepsilon_k}(\tau_1,x) & =1-B_{+}(x)K_{+}(x)\tau_1                             \\
		                                                              & \leq 1-\frac{\Lambda_1-\delta}{\Lambda_1-2\delta} <0,
	\end{align*}
	which contradicts \eqref{Eq:lim D}.
	Therefore, we obtain
	\begin{equation*}
		T^*(\varepsilon)\varepsilon^{p-1} \leq \frac{1}{\Lambda_1-2\delta}.
	\end{equation*}
	By the arbitrariness of $\delta>0$,
	\begin{equation*}
		\limsup_{\varepsilon \rightarrow 0}T^*(\varepsilon)\varepsilon^{p-1} \leq \Lambda_1^{-1}
	\end{equation*}
	holds.
\end{proof}
Next, we estimate the limit inferior of $T^*(\varepsilon)\varepsilon^{p-1}.$
\begin{lem}
	\begin{equation*}
		\liminf_{\varepsilon \rightarrow 0}T^*(\varepsilon)\varepsilon^{p-1} \geq \Lambda_1^{-1}.
	\end{equation*}
\end{lem}
\begin{proof}
	We assume that
	\begin{equation*}
		L \coloneq\liminf_{\varepsilon \rightarrow 0}T^*(\varepsilon)\varepsilon^{p-1} <\Lambda_1^{-1}.
	\end{equation*}
	Then we can take a sequence of positive real numbers $\{\varepsilon_k\}_k$ such that
	\begin{equation*}
		\lim_{k \rightarrow \infty}T^*(\varepsilon_k)\varepsilon_k^{p-1}=L.
	\end{equation*}
	We fix $L<\tau_0<\Lambda_1^{-1}.$
	Then, there exists $k_0 \in \mathbb{N}$ such that
	\begin{equation*}
		T^*(\varepsilon_k)\varepsilon_k^{p-1} < \tau_0
	\end{equation*}
	if $k \geq k_0.$
	That is, $\min(\tau_0,T^*(\varepsilon_k)\varepsilon_k^{p-1})=T^*(\varepsilon_k)\varepsilon_k^{p-1}$ if $k \geq k_0.$
	From the decay assumption $\varphi_x(x), \psi(x) \rightarrow 0$ as $|x| \rightarrow \infty$,
	for any $\eta>0$, there exists $R>0$ such that
	\begin{equation}\label{Eq:F_1 decay}
		\left|\frac{F_j(0)}{\varepsilon_k}\right| \leq \eta
	\end{equation}
	for all $|x| \geq R.$
	We perform the estimates separately for $|x|\geq R$ and $|x|\leq R.$
	By $\gamma(u) \leq C|u|^{p-2}$, \eqref{Eq:u_x} and \eqref{Eq:F_1 decay},
	\begin{align*}
		|D^{(j)}_{\varepsilon_k}(\tau,x)-	1| & =\left|\frac{F_j(0)}{\varepsilon_k} \int_{0}^{\tau}\varepsilon_k^{-(p-2)}\gamma(u(\varepsilon_k^{-(p-1)}\sigma,x_\pm(\varepsilon_k^{-(p-1)}\sigma;0,x)))d\sigma\right| \\
		                                     & \leq C \Lambda_1^{-1}\eta
	\end{align*}
	if $k \geq k_0$ for all $\tau \in [0,T^*(\varepsilon_k)\varepsilon_k^{p-1})$ and $|x| \geq R$.
	Thus,
	\begin{equation}\label{Eq:far}
		D^{(j)}_{\varepsilon_k}(\tau,x)>1-\Lambda_1^{-1}C\eta
	\end{equation}
	for all $(\tau,x) \in [0,T^*(\varepsilon_k)\varepsilon_k^{p-1}) \times \{|x| \geq R\}$.
	Next, we derive the estimate on $[-R,R]$.
	From Lemma 3.3, there exists $k_0 \in \mathbb{N}$ such that
	\begin{equation*}
		|D^{(j)}_{\varepsilon_k}(\tau,x)-D_\pm(\tau,x)|< \eta
	\end{equation*}
	if $k>k_0$ for all $(\tau,x) \in [0,T^*(\varepsilon_k)\varepsilon_k^{p-1}) \times [-R,R]$.
	Therefore,
	\begin{align}
		D^{(j)}_{\varepsilon_k}(\tau,x) & > D_\pm(\tau,x)-\eta    \notag                         \\[4pt]
		                                & \geq (1-\Lambda_1 \tau_0)-\eta \label{Eq: near origin}
	\end{align}
	for all $(\tau,x) \in [0,T^*(\varepsilon_k)\varepsilon_k^{p-1}) \times [-R,R]$.
	In conclusion, by \eqref{Eq:far} and \eqref{Eq: near origin}, we obtain that there exists $k_0 \in \mathbb{N}$ such that
	\begin{equation*}
		D^{(j)}_{\varepsilon_k}(\tau,x) > C
	\end{equation*}
	if $k \geq k_0$ for all $(\tau, x) \in [0,T^*(\varepsilon_k)\varepsilon_k^{p-1}) \times \mathbb{R},$ where $C>0$ is a constant independent of $(\tau, x) \in [0,T^*(\varepsilon_k)\varepsilon_k^{p-1}) \times \mathbb{R}.$
	Thus, $F_j$ is bounded for $(t,x) \in [0,T^*(\varepsilon_k)) \times \mathbb{R}$, which allows the solution to be extended beyond $T^*(\varepsilon_k)$.
	This contradicts the fact that $T^*(\varepsilon_k)$ is the lifespan.
\end{proof}
\begin{proof}[Proof of Theorem 1.1]
	Lemmas 3.4 and 3.5 complete the proof of Theorem 1.1.
\end{proof}

\section{Proof of Theorem 1.2}

\subsection{In the case of $0< \alpha <1$}
\begin{lem}
	We fix $K>0$ arbitrarily.
	There exists $\varepsilon_0>0$ such that if $0 < \varepsilon \leq \varepsilon_0$, then, for $t \in [0,\min (T^*(\varepsilon),T_\varepsilon))$,
	\begin{align}
		\|u(t)\|_{L^\infty}      & \leq C_1\varepsilon \Braket{t}^{1-\alpha},  \label{Eq:u_x a}                       \\
		\|c'(u(t))\|_{L^\infty}  & \leq C_1\varepsilon^{p-2} \Braket{t}^{(1-\alpha)(p-2)}, \label{Eq:c_prime_upper a} \\
		\|c(u(t))-1\|_{L^\infty} & \leq C_1\varepsilon^{p-1} \Braket{t}^{(1-\alpha)(p-1)}, \label{Eq:c_minus_1 a}
	\end{align}
	where $T_\varepsilon=K\varepsilon^{-q}.$
\end{lem}
\begin{proof}
	By the invariance of the Riemann invariants on the characteristics \eqref{Eq:Invariance2},
	\begin{align*}
		2G(u) & =r(0,x_1)-s(0,x_2)                                                                               \\
		      & =\int_{x_2}^{x_1}\varepsilon \psi(y)dy +G(\varepsilon \varphi(x_1))+G(\varepsilon \varphi(x_2)),
	\end{align*}
	where $x_1=x_-(0), x_2=x_+(0),$ that is,
	\begin{align*}
		x=x_1-\int_{0}^{t}c(u(\tau,x_-(\tau)))d\tau, \\
		x=x_2+\int_{0}^{t}c(u(\tau,x_+(\tau)))d\tau.
	\end{align*}
	We use a bootstrap argument for \eqref{Eq:u_x a}.
	Fix $0<T<\min\{T^*(\varepsilon),T_\varepsilon\}$
	and choose $C_1>\|\varphi\|_{L^\infty}$ sufficiently large.
	If \eqref{Eq:u_x a} fails before $T$, let $T_1<T$
	be the first time at which equality holds in \eqref{Eq:u_x a}.
	We derive an improved estimate on $[0,T_1]$
	and obtain a contradiction.
	Regarding the relation between $x_1$ and $x_2$,
	\begin{equation*}
		|x_1-x_2| \leq \int_{0}^{t}|c(u(\sigma,x_-(\sigma)))-1|d\sigma + \int_{0}^{t}|c(u(\sigma,x_+(\sigma)))-1|d\sigma +2t.
	\end{equation*}
	Here we show that
	\begin{equation}
		\label{eq:c-minus-one}
		|c(u)-1|\le C|u|^{p-1}
	\end{equation}
	on \([0,T_1]\times\mathbb R\), provided that
	\(\varepsilon_0>0\) is sufficiently small.
	By the mean value theorem, there exists \(\xi\in(0,1)\)
	such that
	\[
		c(u)=1+c'(\xi u)u.
	\]
	By the assumption \eqref{Eq:unif cvg}, for any \(\eta>0\) and \(M>0\),
	\[
		\rho^{-(p-2)}|c'(\rho\theta)|
		<\eta+a|\theta|^{p-2}
	\]
	for \(|\theta|\le M\), provided that \(\rho>0\)
	is sufficiently small.
	Under the bootstrap assumption,
	\[
		\begin{aligned}
			|u|
			 & \le C_1\varepsilon\langle t\rangle^{1-\alpha} \\
			 & \le C_1\varepsilon
			\langle K\varepsilon^{-q}\rangle^{1-\alpha}      \\
			 & \le C_2\varepsilon^{
				\frac{\alpha}{(1-\alpha)(p-2)+1}}
			\longrightarrow0
			\qquad\text{as }\varepsilon\downarrow0.
		\end{aligned}
	\]
	Thus, taking \(M=1\), \(\rho=|u|\), and
	\(\theta=u/|u|\)  in \eqref{Eq:unif cvg}, we obtain, for \(u\ne0\),
	\[
		|u|^{-(p-2)}
		\left|c'\!\left(|u|\frac{u}{|u|}\right)\right|
		<\eta+a,
	\]
	provided that \(\varepsilon_0>0\) is sufficiently small.
	The same estimate applies to \(\xi u\). Hence,
	\[
		|c(u)-1|
		=|c'(\xi u)|\,|u|
		\le C|\xi u|^{p-2}|u|
		\le C|u|^{p-1}.
	\]
	Since \(c(0)=1\), this also holds for \(u=0\),
	which proves \eqref{eq:c-minus-one}.
	Thus,
	\begin{align*}
		|x_1-x_2| & \leq CC_1^{p-1}\varepsilon^{p-1}\int_{0}^{t}\Braket{\tau}^{(1-\alpha)(p-1)}d\tau+2t \\
		          & \leq C \varepsilon^{p-1}\varepsilon^{-q(1-\alpha)(p-1)}t + 2t                       \\
		          & = C \varepsilon^{(p-1)(1-q(1-\alpha))}t + 2t                                        \\
		          & = C \varepsilon^{\frac{(p-1)\alpha}{(1-\alpha)(p-2)+1}}t +2t.
	\end{align*}
	By taking $\varepsilon_0>0$ sufficiently small,
	\begin{equation*}
		|x_1-x_2| \leq 3t.
	\end{equation*}
	We estimate each term of $2G(u)$.
	\begin{align*}
		|G(\varepsilon \varphi(x_1))+G(\varepsilon \varphi(x_2))| & \leq 2\|G(\varepsilon \varphi )\|_{L^\infty} \leq C \varepsilon, \\[7pt]
		\left|\int_{x_2}^{x_1} \varepsilon \psi_1(y)dy\right| \leq \varepsilon \|\psi_1\|_{L^1}.
	\end{align*}
	Hence,
	\begin{align*}
		|2G(u)| & \leq C\varepsilon + \varepsilon\int_{x_2}^{x_1}|\psi_2(y)|dy              \\
		        & \leq C\varepsilon + \varepsilon \int_{x_2}^{x_1} \Braket{y}^{-\alpha}dy   \\
		        & \leq C \varepsilon\left(1+\int_{0}^{\frac{3}{2}t}(1+y)^{-\alpha}dy\right) \\
		        & \leq C \varepsilon(1+ \Braket{t}^{1-\alpha})                              \\[8pt]
		        & \leq C \varepsilon \Braket{t}^{1-\alpha}.
	\end{align*}
	Moreover, by $G(u)=u+\mathcal{O}(|u|^p)$,
	\begin{equation*}
		\|u(t)\|_{L^\infty} \leq C\varepsilon \Braket{t}^{1-\alpha}.
	\end{equation*}
	By the definition of $T_1$, we have
	\begin{equation*}
		\|u(t)\|_{L^\infty} \leq C_1\varepsilon \Braket{t}^{1-\alpha},
	\end{equation*}
	and
	\begin{equation*}
		\|u(T_1)\|_{L^\infty} = C_1\varepsilon \Braket{T_1}^{1-\alpha}.
	\end{equation*}
	Therefore,
	\begin{equation*}
		C_1 \leq C
	\end{equation*}
	holds. On the other hand, choosing $C_1 > \max(\|\varphi\|_{L^\infty}, C)$
	in advance yields $C_1 > C$, which is a contradiction.
	Therefore, \eqref{Eq:u_x a} holds for any $t \in [0,\min(T^*(\varepsilon), T_\varepsilon)).$
	It is clear that \eqref{Eq:c_prime_upper a} and \eqref{Eq:c_minus_1 a} hold.
\end{proof}

We follow the argument in the proof of Theorem 1.1.
We set
\begin{equation*}
	\tilde{D}^{(j)}_\varepsilon(\tau,x)=D^{(j)}(\varepsilon^{-q}\tau,x).
\end{equation*}
\begin{lem}
	We fix $R>0, \tau_0 >0$ and set $I_\varepsilon=[0,\min(\tau_0, T^*(\varepsilon)\varepsilon^q)).$
	Then
	\begin{align*}
		\sup_{(\tau,x)\in I_\varepsilon \times [-R,R]}|\tilde{D}^{(1)}_\varepsilon(\tau,x) - \tilde{D}_{+}(\tau,x)| \rightarrow 0, \\
		\sup_{(\tau,x)\in I_\varepsilon \times [-R,R]}|\tilde{D}^{(2)}_\varepsilon(\tau,x) - \tilde{D}_{-}(\tau,x)| \rightarrow 0,
	\end{align*}
	as $\varepsilon \rightarrow 0$, where
	\begin{equation*}
		\tilde{D}_{\pm}(\tau,x)=1-B_\pm(x)J(\tau), \hspace{10pt} J(\tau)=\frac{2^{-\alpha(p-2)}a}{2k(1-\alpha)^{p-2}}\tau^k.
	\end{equation*}
\end{lem}
\begin{proof}
	We show that, for any $\eta>0$, there exists $\varepsilon_0>0$ such that
	\begin{align*}
		|\tilde{D}^{(1)}_\varepsilon(\tau,x) - \tilde{D}_{+}(\tau,x)| < \eta, \\
		|\tilde{D}^{(2)}_\varepsilon(\tau,x) - \tilde{D}_{-}(\tau,x)| < \eta,
	\end{align*}
	for $(\tau,x) \in [0,\min(\tau_0, T^*(\varepsilon)\varepsilon^q)) \times [-R,R].$
	We fix $\delta>0$ arbitrarily and set $I_{\delta,\varepsilon} = [\delta, \min(\tau_0, T^*(\varepsilon)\varepsilon^q))$.
	We assume that $(\tau, x) \in I_{\delta,\varepsilon} \times [-R,R]$ and set $t=\varepsilon^{-q}\tau.$
	We consider $x_-(\sigma)$ through $(0,x)$ and
	$x_+(\sigma)$ intersecting $x_-(\sigma)$ at time $t$, and let $(t,X)$ be their point of intersection.
	We set $x_+(0)=y_\varepsilon$, since the initial value of $x_+(\sigma)$ depends on $\varepsilon$.
	From $2G(u)=r-s$ and \eqref{Eq:Invariance2},
	\begin{align*}
		2G(u(t,X)) & =r(0,x)-s(0,y_\varepsilon)                                                                                        \\
		           & = \int_{y_\varepsilon}^{x} \varepsilon \psi(y)dy+G(\varepsilon \varphi(x))+G(\varepsilon \varphi(y_\varepsilon)).
	\end{align*}
	Then
	\begin{equation*}
		y_\varepsilon \leq R - \varepsilon^{-q}\delta.
	\end{equation*}
	First, we consider $\varepsilon^{-1+q(1-\alpha)}G(u).$
	It is clear that there exists $\varepsilon_0>0$ such that
	\begin{align*}
		\varepsilon^{-1+q(1-\alpha)}\max(|G(\varepsilon \varphi (y_\varepsilon))|,|G(\varepsilon \varphi (x))|) < \eta
	\end{align*}
	for all $|x| \leq R.$
	For the first term of $G(u)$, we have
	\begin{align*}
		\int_{y_\varepsilon}^{x}\psi(y)dy & =\int_{y_\varepsilon}^{x}\psi_1(y)dy+ \int_{y_\varepsilon}^{x}\psi_2(y)dy \\
		                                  & =E_\varepsilon(x)+\int_{y_{\varepsilon}}^{0}\Braket{y}^{-\alpha}dy,
	\end{align*}
	where
	\begin{align*}
		|\varepsilon^{q(1-\alpha)}E_\varepsilon(x)| & =\left|\varepsilon^{q(1-\alpha)} \left(\int_{y_\varepsilon}^{x}\psi_1(y)dy+\int_{0}^{x}\psi_2(y)dy\right)\right| \\
		                                            & \leq \varepsilon^{q(1-\alpha)}\left(\|\psi_1\|_{L^1}+\int_{0}^{|x|}\Braket{y}^{-\alpha} dy\right)                \\
		                                            & \leq \varepsilon^{q(1-\alpha)}\left(\|\psi_1\|_{L^1}+\int_{0}^{R}\Braket{y}^{-\alpha} dy\right).
	\end{align*}
	We decompose
	\begin{equation*}
		\int_{y_\varepsilon}^{0}\Braket{y}^{-\alpha}dy=(-y_\varepsilon)^{1-\alpha} \left( (-y_\varepsilon)^{-(1-\alpha)}\int_{y_\varepsilon}^{0}\Braket{y}^{-\alpha}dy \right) \eqqcolon (-y_\varepsilon)^{1-\alpha}  Y_\varepsilon.
	\end{equation*}
	Then, applying l'H\^{o}pital's theorem yields $Y_\varepsilon \rightarrow \frac{1}{1-\alpha}.$
	Furthermore,
	\begin{align*}
		y_\varepsilon & = x -\int_{0}^{t}c(u(\sigma,x_-(\sigma)))d\sigma - \int_{0}^{t}c(u(\sigma,x_+(\sigma)))d\sigma               \\
		              & = x - 2t -\int_{0}^{t}(c(u(\sigma,x_-(\sigma)))-1)d\sigma - \int_{0}^{t}(c(u(\sigma,x_+(\sigma)))-1)d\sigma,
	\end{align*}
	and
	\begin{equation*}
		\varepsilon^q (-y_\varepsilon) = -\varepsilon^q x +2\tau + \varepsilon^q \left(\int_{0}^{t}(c(u(\sigma,x_-(\sigma)))-1)d\sigma + \int_{0}^{t}(c(u(\sigma,x_+(\sigma)))-1)d\sigma\right).
	\end{equation*}
	For the integral term,
	\begin{align*}
		\left| \int_{0}^{t}(c(u(\sigma,x_\pm(\sigma)))-1)d\sigma \right| & \leq C \varepsilon^{p-1} \int_{0}^{t}\Braket{\sigma}^{(p-1)(1-\alpha)}d\sigma \\
		                                                                 & \leq C \varepsilon^{p-1} t^{(1-\alpha)(p-1)+1}                                \\[6pt]
		                                                                 & \leq C \varepsilon^{p-1} \varepsilon^{-q(1-\alpha)(p-1)-q}.
	\end{align*}
	Thus, we have
	\begin{equation*}
		\varepsilon^q (-y_\varepsilon) \rightarrow 2\tau
	\end{equation*}
	from $q+(p-1)-q(1-\alpha)(p-1)-q=\alpha q>0$, and
	\begin{equation*}
		\varepsilon^{q(1-\alpha)}(-y_\varepsilon)^{1-\alpha} = (\varepsilon^q (-y_\varepsilon))^{1-\alpha} \rightarrow (2\tau)^{1-\alpha}.
	\end{equation*}
	Therefore, there exists $\varepsilon_0>0$ such that
	\begin{equation*}
		\left|\varepsilon^{-1+q(1-\alpha)}G(u(t,X)) - \frac{2^{-\alpha}}{1-\alpha}\tau^{1-\alpha}\right| < \eta,
	\end{equation*}
	and
	\begin{equation*}
		\left|\varepsilon^{-1+q(1-\alpha)}u(t,X) - \frac{2^{-\alpha}}{1-\alpha}\tau^{1-\alpha}\right| < \eta,
	\end{equation*}
	if $0<\varepsilon<\varepsilon_0$ for all $(\tau,x) \in I_{\delta,\varepsilon} \times [-R,R],$
	by $G(u)=u+ \mathcal{O}(\varepsilon^{p(1-q(1-\alpha))})$.
	Similarly to the proof of Theorem 1.1, we have
	\begin{align*}
		\left|\varepsilon^{(-1+q(1-\alpha))(p-2)}\gamma(u(t,X))-\Gamma(\tau)\right| < \eta,
	\end{align*}
	where $\Gamma(\tau)=\frac{2^{-\alpha(p-2)}a}{2(1-\alpha)^{p-2}}\tau^{(1-\alpha)(p-2)}.$
	For the exponent, noting that
	\begin{equation*}
		1-q+(1-q(1-\alpha))(p-2)=0,
	\end{equation*}
	we have
	\begin{align*}
		|\tilde{D}^{(1)}_\varepsilon(\tau,x)-\tilde{D}_{+}(\tau,x)| & = \left|\frac{F_1(0)}{\varepsilon}\int_{0}^{\tau}\varepsilon^{(-1+q(1-\alpha))(p-2)}\gamma(u)d\sigma-B_+(x)J(\tau)\right|                    \\[7pt]
		                                                            & \leq \left| \frac{F_1(0)}{\varepsilon} \right|\int_{0}^{\tau}\left|\varepsilon^{(-1+q(1-\alpha))(p-2)}\gamma(u)-\Gamma(\sigma)\right|d\sigma \\
		                                                            & +\left| \frac{F_1(0)}{\varepsilon}-B_+(x) \right|J(\tau).
	\end{align*}
	For any \(\eta>0\), by choosing \(\delta>0\)
	sufficiently small, we have
	\begin{equation*}
		\int_0^\delta
		\left|
		\varepsilon^{(-1+q(1-\alpha))(p-2)}\gamma(u)
		-\Gamma(\sigma)
		\right|\,d\sigma<\eta
	\end{equation*}
	by the uniform boundedness of
	\(\varepsilon^{(-1+q(1-\alpha))(p-2)}\gamma(u)\)
	and \(\Gamma(\sigma)\).
	For this fixed \(\delta\), we also have
	\begin{equation*}
		\int_\delta^\tau
		\left|
		\varepsilon^{(-1+q(1-\alpha))(p-2)}\gamma(u)
		-\Gamma(\sigma)
		\right|\,d\sigma<\eta
	\end{equation*}
	by the uniform convergence above, provided that
	\(\varepsilon_0>0\) is sufficiently small.
	The estimates for the remaining terms are similar
	to those in Theorem 1.1.
	Therefore,
	\begin{equation*}
		|\tilde{D}^{(1)}_\varepsilon(\tau,x)-\tilde{D}_{+}(\tau,x)| < \eta,
	\end{equation*}
	if $0<\varepsilon< \varepsilon_0$ uniformly with respect to $(\tau,x) \in I_\varepsilon \times [-R, R].$
	The same holds for $\tilde{D}^{(2)}_\varepsilon(\tau,x).$
\end{proof}
\begin{proof}[Proof of Theorem 1.2 in the case of \,\,$0<\alpha <1$]
	The upper and lower lifespan bounds follow from
	the arguments of Lemmas 3.4 and 3.5, respectively,
	using Lemmas 4.1 and 4.2 in place of
	Proposition 3.2 and Lemma 3.3.
\end{proof}

\subsection{In the case of $\alpha=1$}
\begin{lem}
	We fix $K>0$ arbitrarily.
	There exists $\varepsilon_0>0$ such that if $0 < \varepsilon \leq \varepsilon_0$, then, for $t \in [0,\min (T^*(\varepsilon),T_\varepsilon))$,
	\begin{align}
		\|u(t)\|_{L^\infty}      & \leq C_1\varepsilon \log(2+t),  \label{Eq:u_x log}                        \\
		\|c'(u(t))\|_{L^\infty}  & \leq C_1\varepsilon^{p-2} (\log(2+t))^{p-2}, \label{Eq:c_prime_upper log} \\
		\|c(u(t))-1\|_{L^\infty} & \leq C_1\varepsilon^{p-1} (\log(2+t))^{p-1}, \label{Eq:c_minus_1 log}
	\end{align}
	where $T_\varepsilon=K\varepsilon^{-(p-1)}|\log \varepsilon|^{-(p-2)}$.
\end{lem}
\begin{proof}
	This is completely analogous to the case $0<\alpha<1$.
	In the estimation of $G(u)$, a logarithmic term arises since
	\begin{equation*}
		|2G(u)| \leq C\varepsilon + \varepsilon \int_{x_2}^{x_1}\Braket{y}^{-1}dy.
	\end{equation*}

\end{proof}
As in the case of $0<\alpha<1$, we employ a similar argument
to that in the proof of Theorem 1.1.
We set
\begin{equation*}
	\mathcal{D}^{(j)}_\varepsilon(\tau,x)= D^{(j)}(\varepsilon^{-(p-1)}|\log\varepsilon|^{-(p-2)}\tau,x).
\end{equation*}
\begin{lem}
	We fix $R>0, \tau_0 >0$ and set $I_\varepsilon=[0,\min(\tau_0, T^*(\varepsilon)\varepsilon^{p-1}|\log \varepsilon|^{p-2})).$
	Then
	\begin{align*}
		\sup_{(\tau,x)\in I_\varepsilon \times [-R,R]}|\mathcal{D}^{(1)}_\varepsilon(\tau,x) - \mathcal{D}_{+}(\tau,x)| \rightarrow 0, \\
		\sup_{(\tau,x)\in I_\varepsilon \times [-R,R]}|\mathcal{D}^{(2)}_\varepsilon(\tau,x) - \mathcal{D}_{-}(\tau,x)| \rightarrow 0,
	\end{align*}
	as $\varepsilon \rightarrow 0$, where
	\begin{equation*}
		\mathcal{D}_{\pm}(\tau,x)=1-B_\pm(x)L\tau, \hspace{10pt} L=\frac{a(p-1)^{p-2}}{2^{p-1}}.
	\end{equation*}
\end{lem}
\begin{proof}
	The proof is similar to those of Lemmas 3.3 and 4.2.
	We fix $\delta, \eta>0$ arbitrarily and set $I_{\delta,\varepsilon}=[\delta, \min(\tau_0, T^*(\varepsilon)\varepsilon^{p-1}|\log \varepsilon|^{p-2})).$
	We assume that $(\tau,x) \in I_{\delta,\varepsilon} \times [-R,R]$  and set $t=\varepsilon^{-(p-1)}|\log \varepsilon|^{-(p-2)}\tau$.
	We consider $x_-(\sigma)$ through $(0,x)$ and
	$x_+(\sigma)$ intersecting $x_-(\sigma)$ at time $t$, and let $(t,X)$ be their point of intersection.
	We set $x_+(0)=y_\varepsilon$.
	By the invariance of Riemann invariants on the characteristics \eqref{Eq:Invariance2},
	\begin{align*}
		2G(u) & =r(0,x)-s(0,y_\varepsilon)                                                                                        \\
		      & =\int_{y_\varepsilon}^{x}\varepsilon \psi(y)dy +G(\varepsilon \varphi(x))+ G(\varepsilon \varphi(y_\varepsilon)).
	\end{align*}
	Then
	\begin{equation*}
		y_\varepsilon \leq x-\varepsilon^{-(p-1)}|\log \varepsilon|^{-(p-2)}\tau.
	\end{equation*}
	First, we consider $\varepsilon^{-1}|\log \varepsilon|^{-1}G(u).$
	It is clear that there exists $\varepsilon_0>0$ such that
	\begin{equation*}
		\varepsilon^{-1}|\log \varepsilon|^{-1}\max(|G(\varepsilon \varphi (y_\varepsilon))|,|G(\varepsilon \varphi (x)) |) < \eta.
	\end{equation*}
	For the first term of $G(u)$, we have
	\begin{align*}
		\int_{y_\varepsilon}^{x}\psi(y)dy & =\int_{y_\varepsilon}^{x}\psi_1(y)dy+ \int_{y_\varepsilon}^{x}\psi_2(y)dy \\
		                                  & =E_\varepsilon(x)+\int_{y_{\varepsilon}}^{0}\Braket{y}^{-1}dy,
	\end{align*}
	where
	\begin{align*}
		||\log \varepsilon|^{-1} E_\varepsilon(x)| & =\left||\log \varepsilon|^{-1} \left(\int_{y_\varepsilon}^{x}\psi_1(y)dy+\int_{0}^{x}\psi_2(y)dy\right)\right| \\
		                                           & \leq |\log \varepsilon|^{-1}\left(\|\psi_1\|_{L^1}+\int_{0}^{R}\Braket{y}^{-1} dy\right).
	\end{align*}
	We decompose
	\begin{equation*}
		\int_{y_\varepsilon}^{0}\Braket{y}^{-1}dy=\log (-y_\varepsilon)\left((\log (-y_\varepsilon))^{-1}\int_{y_\varepsilon}^{0}\Braket{y}^{-1}dy \right) \eqqcolon \log(-y_\varepsilon) Y_\varepsilon.
	\end{equation*}
	Then, applying l'H\^{o}pital's theorem yields $Y_\varepsilon \rightarrow 1.$
	Furthermore,
	\begin{align*}
		y_\varepsilon & = x -\int_{0}^{t}c(u(\sigma,x_-(\sigma)))d\sigma - \int_{0}^{t}c(u(\sigma,x_+(\sigma)))d\sigma                                                                 \\
		              & = x - 2t -\int_{0}^{t}(c(u(\sigma,x_-(\sigma)))-1)d\sigma - \int_{0}^{t}(c(u(\sigma,x_+(\sigma)))-1)d\sigma                                                    \\
		              & = x - 2\varepsilon^{-(p-1)}|\log \varepsilon|^{-(p-2)}\tau - \int_{0}^{t}(c(u(\sigma,x_-(\sigma)))-1)d\sigma - \int_{0}^{t}(c(u(\sigma,x_+(\sigma)))-1)d\sigma \\
		              & = - 2\varepsilon^{-(p-1)}|\log \varepsilon|^{-(p-2)}\tau \cdot R_\varepsilon(t,x),
	\end{align*}
	where
	\begin{equation*}
		R_\varepsilon(t,x)=1-\frac{x-\int_{0}^{t}(c(u(\sigma,x_-(\sigma)))-1)d\sigma - \int_{0}^{t}(c(u(\sigma,x_+(\sigma)))-1)d\sigma }{2\varepsilon^{-(p-1)}|\log \varepsilon|^{-(p-2)}\tau}.
	\end{equation*}
	Thus,
	\begin{align*}
		|\log \varepsilon|^{-1} \log(-y_\varepsilon) & =|\log \varepsilon|^{-1} \log(2\tau)+(p-1)                                                              \\
		                                             & -(p-2)|\log \varepsilon|^{-1} \log |\log \varepsilon| + |\log\varepsilon|^{-1} \log R_\varepsilon(t,x).
	\end{align*}
	We set $u_\pm = u(\sigma, x_\pm(\sigma)).$
	For the remainder term $R_\varepsilon(t, x)$,
	since \begin{equation*}
		|R_\varepsilon(t,x) - 1| \leq \frac{|x| + \sum\limits_{\pm}\int_{0}^{t}|c(u_\pm)-1|d\sigma}{2 \varepsilon^{-(p-1)}|\log \varepsilon|^{-(p-2)}\delta}
	\end{equation*}
	and
	\begin{align*}
		\int_{0}^{t}|c(u_\pm) - 1|d\sigma & \leq C \varepsilon^{p-1}|\log \varepsilon|^{p-1}t                                                            \\
		                                  & \leq C \varepsilon^{p-1}|\log \varepsilon|^{p-1} \cdot \varepsilon^{-(p-1)}|\log \varepsilon|^{-(p-2)}\tau_0 \\[4pt]
		                                  & = C |\log \varepsilon|\tau_0,
	\end{align*}
	we have
	\begin{equation*}
		R_\varepsilon(t,x) \rightarrow 1.
	\end{equation*}
	Consequently,
	\begin{equation*}
		\left|\varepsilon^{-1}|\log \varepsilon|^{-1} G(u(t,X)) - \frac{p-1}{2}\right| < \eta,
	\end{equation*}
	and
	\begin{equation*}
		\left|\varepsilon^{-1}|\log \varepsilon|^{-1} u(t,X) - \frac{p-1}{2}\right|<\eta.
	\end{equation*}
	Therefore,
	\begin{align*}
		|\varepsilon^{-(p-2)}|\log \varepsilon|^{-(p-2)}\gamma(u(t,X)) - L| < \eta,
	\end{align*}
	and
	\begin{equation*}
		\sup_{(\tau,x)\in I_\varepsilon \times [-R,R]}|\mathcal{D}^{(1)}_\varepsilon(\tau,x) - \mathcal{D}_{+}(\tau,x)| < \eta
	\end{equation*}
	can be shown by the same argument as in the case $0 < \alpha < 1$.
	The same holds for $\mathcal{D}^{(2)}_\varepsilon(\tau,x).$
\end{proof}
\begin{proof}[Proof of Theorem 1.2 in the case of \,\,$\alpha=1$]
	We argue as in the case \(0<\alpha<1\),
	using Lemmas 4.3 and 4.4 in place of
	Lemmas 4.1 and 4.2.
\end{proof}

Yuusuke Sugiyama

Department of Mathematics,

Tokyo University of Science,

Kagurazaka 1-3, Shinjuku-ku,

Tokyo 162-8601, Japan

E-mail: sugiyama.y@rs.tus.ac.jp

\vspace{20pt}
Taro Yamanoi

Department of Mathematics,

Faculty of Science,

Graduate School of Tokyo University of Science,

Kagurazaka 1-3, Shinjuku-ku,

Tokyo 162-8601, Japan

E-mail: 1126709@ed.tus.ac.jp

\end{document}